\documentclass[12pt]{article}
\usepackage{mathrsfs}

\usepackage{amsmath,amssymb,amsfonts,theorem,makeidx,latexsym,epsfig,color,subfigure,algorithm,algorithmic}\usepackage[numbers,sort&compress]{natbib}

\newtheorem{defn}{Definition}[section]

\newtheorem{lem}[defn]{Lemma}

{\theorembodyfont{\rmfamily}

}

\newtheorem{thm}[defn]{Theorem}

\newtheorem{obs}[defn]{Observation}

\newtheorem{prop}[defn]{Proposition}

\newtheorem{cor}[defn]{Corollary}

\numberwithin{equation}{section}

\def\bp{{\noindent\bf Proof. \ }}

\def\ep{\hfill$\square$\par\bigskip}

\title{Double domination in maximal outerplanar graphs \footnote{The research is supported by NSFC (No. 11301440),
Natural Science Foundation of Fujian Province (CN)(2015J05017)}}

\author{{\large Wei Zhuang\thanks{Corresponding author; E-mail: zhuangweixmu@163.com}} \\
{\it \normalsize School of Applied Mathematics},
\\{\it \normalsize  Xiamen University of Technology, Xiamen 361024, P.R.China }}

\date{}

\begin{document}

\maketitle

\begin{abstract} In a graph $G$, a vertex dominates itself and its neighbors.
A subset $S\subseteq V(G)$ is said to be a double dominating set of
$G$ if $S$ dominates every vertex of $G$ at least twice. The double
domination number $\gamma_{\times 2}(G)$ is the minimum cardinality
of a double dominating set of $G$. We show that if $G$ is a maximal
outerplanar graph on $n\geq 3$ vertices, then $\gamma_{\times
2}(G)\leq \lfloor \frac{2n}{3}\rfloor$. Further, if $n\geq 4$, then
$\gamma_{\times 2}(G)\leq \min \{\lfloor \frac{n+t}{2}\rfloor,
n-t\}$, where $t$ is the number of vertices of degree $2$ in $G$.
These bounds are shown to be tight. In addition, we also study the
case that $G$ is a striped maximal outerplanar graph.

\end{abstract}

\begin{minipage}{150mm}

\end{minipage}

\section{ Introduction}
For a simple graph $G=(V, E)$, $V$ and $E$ are the sets of vertices
and edges of $G$, respectively. We denote by $|V(G)|$ the
\emph{order} of a graph $G$. If the graph $G$ is clear from the
context, we simply write $|G|$ or $n$ rather than $|V(G)|$. For a
vertex $v\in V(G)$, let $N_G(v)$ and $N_G[v]$ denote the \emph{open
neighborhood} and the \emph{closed neighborhood} of $v$,
respectively; thus $N_G(v)=\{u \mid uv\in E(G)\}$ and
$N_G[v]=\{v\}\cup N_G(v)$. A graph $G$ is \emph{outerplanar} if it
has embedding in the plane such that all vertices belong to the
boundary of its outer face. An outerplanar graph $G$ is
\emph{maximal} if $G+uv$ is not outerplanar for any two nonadjacent
vertices $u$ and $v$.

A \emph{dominating set} of a graph $G$ is a set $S\subseteq V(G)$
such that every vertex in $G$ is either in $S$ or is adjacent to a
vertex in $S$. A set $D\subseteq V(G)$ is a double dominating set of
$G$ if every vertex in $V(G)-D$ has at least two neighbors in $D$
and every vertex of $D$ has a neighbor in $D$. The \emph{domination
number} (\emph{double domination number}, respectively) of $G$,
denoted by $\gamma(G)$ ($\gamma_{\times 2}(G)$, respectively), is
the minimum cardinality of a dominating set (double dominating set,
respectively). A double dominating set of $G$ of cardinality
$\gamma_{\times 2}(G)$ is called a $\gamma_{\times 2}(G)$-set. If
the graph $G$ is clear from the context, we simply write
$\gamma_{\times 2}$-set rather than $\gamma_{\times 2}(G)$-set. We
say a vertex $v$ in $G$ is double dominated, by a set $S$, if
$|N_G[v]\cap S|\geq 2$.

The concept of double domination was introduced by Harary and Haynes
\cite{Harary} and further studied in, for example, \cite{Blidia,
Chellali, Chellali1, Hajian, Harant, Kaemawichanurat, Mart}.

 In 1996,
Matheson and Tarjan \cite{Matheson} proved that any triangulated
disc $G$ with $n$ vertices satisfies $\gamma(G)\leq \lfloor
\frac{n}{3}\rfloor$, and conjectured that $\gamma(G)\leq \lfloor
\frac{n}{4}\rfloor$ for every $n$-vertex triangulation $G$ with
sufficiently large $n$. In 2013, Campos and Wakabayashi investigated
this problem for maximal outerplanar graphs and showed in
\cite{Campos} that if $G$ is a maximal outerplanar graph of order
$n$, then $\gamma(G)\leq \frac{n+k}{4}$ where $k$ is the number of
vertices of degree $2$ in $G$. Tokunaga proved the same result
independently in \cite{Tokunaga}. Li et al.\cite{Li} improved the
result by showing that $\gamma(G)\leq \frac{n+t}{4}$, where $t$ is
the number of pairs of consecutive $2$-degree vertices with distance
at least $3$ on the outer cycle. For results on other types of
domination in maximal outerplanar graphs, we refer the reader to
\cite{Araki, Borg, Dorfling1, Dorfling2, Henning, Tokunaga1, Zhu,
Lem}.

In this paper, we apply the idea of coloring to investigated the
same question for double domination. We show that if $G$ is a
maximal outerplanar graph on $n\geq 3$ vertices, then
$\gamma_{\times 2}(G)\leq \lfloor \frac{2n}{3}\rfloor$. Further, if
$n\geq 4$, then $\gamma_{\times 2}(G)\leq \min \{\lfloor
\frac{n+t}{2}\rfloor, n-t\}$, where $t$ is the number of vertices of
degree $2$ in $G$. These bounds are shown to be tight. In addition,
we also study the case that $G$ is a striped maximal outerplanar
graph.

\section{Preliminaries}
A maximal outerplanar graph $G$ can be embedded in the plane such
that the boundary of the outer face is a Hamiltonian cycle and each
inner face is a triangle. A maximal outerplanar graph embedded in
the plane is called a \emph{maximal outerplane graph}. For such an
embedding of $G$, we denote by $H_G$ the Hamiltonian cycle which is
boundary of the outer face. An inner face of a maximal outerplane
graph $G$ is an \emph{internal triangle} if it is not adjacent to
outer face. A maximal outerplane graph without internal triangles is
called \emph{stripped}. The following results are useful for our
study.

\begin{prop}[\cite{Campos}]
Let $G$ be a maximal outerplanar graph of order $n\geq 4$. If $G$
has $k$ internal triangles, then $G$ has $k+2$ vertices of degree
$2$.
\end{prop}

\begin{prop}[\cite{Dorfling1}]
If $G$ is a maximal outerplanar graph of order $n\geq 4$, then the
set of vertices of $G$ of degree $2$ is an independent set of $G$ of
size at most $\frac{n}{2}$.
\end{prop}

\begin{prop}[\cite{Harary}]
Let $C$ be a cycle of order $n\geq 3$, then $\gamma_{\times
2}(C)=\lceil \frac{2n}{3} \rceil$.
\end{prop}

\section{Main result}

Let $G$ be any maximal outerplanar graph. We know that there exists
a Hamiltonian cycle in $G$, say $H_G$. It is easy to see that a
$\gamma_{\times 2}$-set of $H_G$ is a double dominating set of $G$.
Hence, the following result is immediate from Proposition~2.3.

\begin{obs}
Let $G$ be a maximal outerplanar graph of order $n\geq 3$, then
$\gamma_{\times 2}(G)\leq \lceil \frac{2n}{3}\rceil$.
\end{obs}

This conclusion is obvious, so we are ready to improve it. Next, we
will give some definitions which are helpful for our investigations.

A \emph{minor} of a graph $G$ is a graph which can be obtained from
$G$ by deleting vertices and deleting or contracting edges. Given a
graph $H$, a graph $G$ is \emph{$H$-minor free} if no minor of $G$
is isomorphic to $H$. A $K_4$-minor free graph $G$ is \emph{maximal}
if $G+uv$ is not a $K_4$-minor free graph for any two non-adjacent
vertices $u$ and $v$ of $G$. A \emph{$2$-tree} is defined
recursively as follows. A single edge is a 2-tree. Any graph
obtained from a 2-tree by adding a new vertex and making it adjacent
to the end vertices of an existing edge is also a 2-tree. It is well
known that the maximal $K_4$-minor free graphs are exactly the
2-trees. We know that a maximal outerplanar graph must be a maximal
$K_4$-minor free graph. Hence, before giving the upper bound for the
double domination number of a maximal outerplanar graph, we intend
to prove the stronger result as follows.

\begin{thm}
Let $G$ be a maximal $K_4$-minor free graph of order $n\geq 3$, then
$\gamma_{\times 2}(G)\leq \lfloor \frac{2n}{3}\rfloor$.
\end{thm}

\bp Let $G_i=G_{i-1}-v_i$ for $i=1, 2, \cdots, n-3$, where $v_i$ is
a vertex of degree 2 in $G_{i-1}$, $G_0=G$. Since a maximal
$K_4$-minor free graph is exactly a 2-tree and according to the
definition of 2-tree, we have that $G_{n-3}$ is a $K_3$. We can give
a proper 3-coloring to $G_0$ by assigning color 1, 2 and 3 to the
three vertices of $G_{n-3}$ respectively and color $v_i$ with the
remaining color which does not appear in $N_{G_{i-1}}(v_i)$ at each
stage. Let $V_t$ be the set of vertices assigned color $t$, where
$t=1, 2, 3$. Choosing two suitable color classes, say $V_1$ and
$V_2$, such that $|V_1|+|V_2|\leq \lfloor \frac{2}{3}|G_0|\rfloor$.
Let $S=V_1\cup V_2$. Next, we will show that $S$ is a double
dominating set of $G$.

Clearly, each of the three vertices is double dominated by $S\cap
V(G_{n-3})$ in $G_{n-3}$. And in $V(G_{i-1})$ ($i=n-3, n-4, \cdots,
2, 1$), $|N_{G_{i-1}}[v_i]\cap (S\cap V(G_{i-1}))|=2$. It means that
each vertex of $G$ is double dominated by $S$. \ep

Based on the above theorem, we have that the following conclusion.

\begin{cor}
Let $G$ be a maximal outerplanar graph of order $n\geq 3$, then
$\gamma_{\times 2}(G)\leq \lfloor \frac{2n}{3}\rfloor$.
\end{cor}

In order to show that this bound is tight. We are ready to construct
an infinite family of graphs as follows. Let $H$ be any maximal
outerplanar graph, and $C=a_1a_2\cdots a_{2k-1}a_{2k}a_1$ is the
unique Hamiltonian cycle of $H$. We note that $C$ is the boundary of
the outer face of $H$. Let $G_H$ be the graph obtained from $H$ by
adding $k$ new vertices $u_1, u_2, \cdots, u_k$ and $2k$ new edges
$u_1a_1, u_1a_2, u_2a_3, u_2a_4, \cdots, u_ka_{2k-1}, u_ka_{2k}$
(see Fig.1, the non-triangular face can be triangulated in any way
so as to obtain a maximal outerplanar graph). Let $\mathscr{U}$ be a
family consisting of all such graph $G_H$. Take a $\gamma_{\times
2}$-set of $G_H$, say $S$. Let $v$ be a vertex of degree two in
$G_H$, it is easy to see that $|N_{G_H}[v]\cap S|\geq 2$. It means
that $|S|\geq 2k=\frac{2|G_H|}{3}$. On the other hand, it follows
from Corollary~3.3 that $|S|\leq \frac{2|G_H|}{3}$. Hence,
$\gamma_{\times 2}(G_H)= \frac{2|G_H|}{3}$.

\begin{center}
  \includegraphics[width=3.5in]{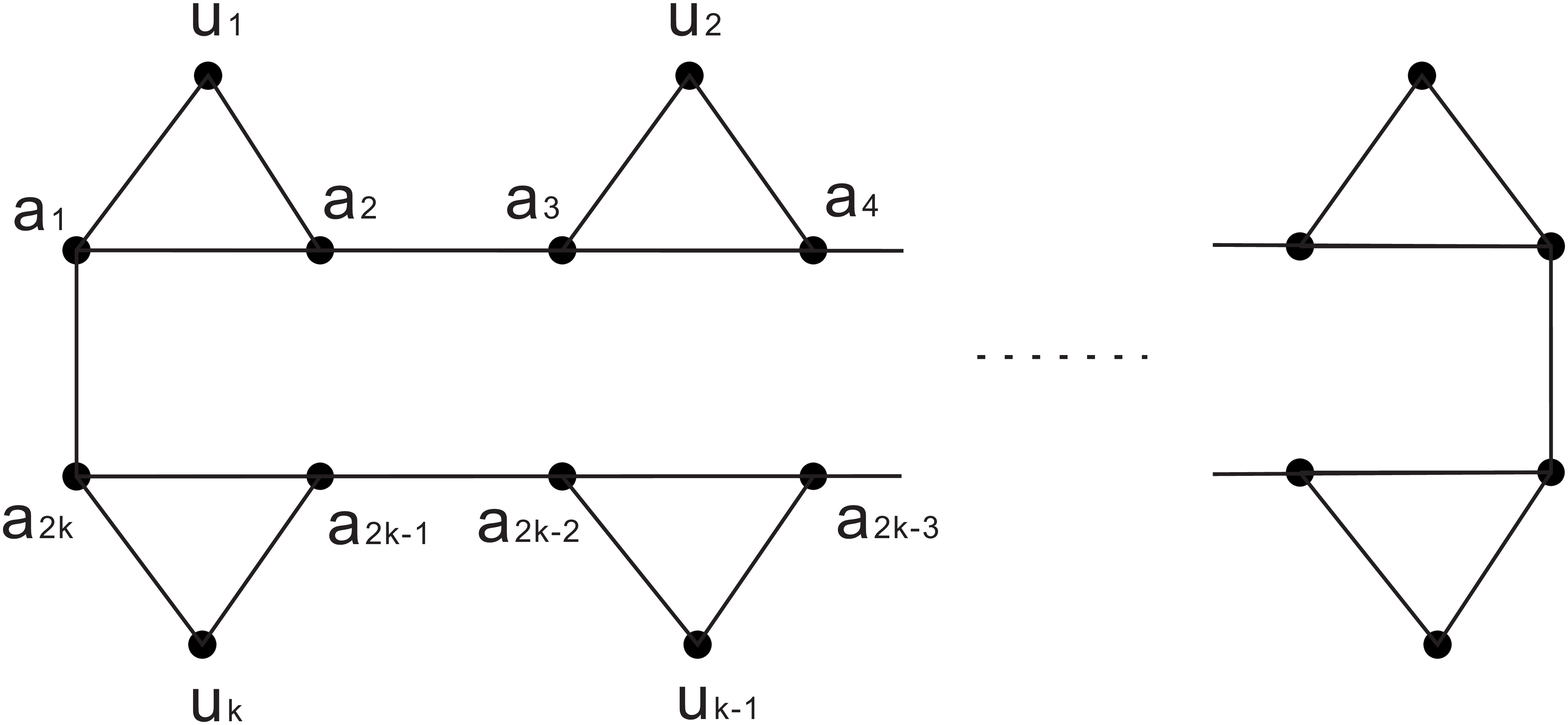}
 \end{center}
\qquad \qquad \qquad \qquad \qquad \qquad \qquad \qquad \qquad
{\small Fig.1} $\\$

But the result of Corollary~3.3 is still not good enough, it only
slightly improves Observation~3.1, so we hope to further improve the
result of Corollary~3.3. Next, we will give two upper bounds in
terms of the order and the number of $2$-degree vertices, and the
first one is inspired by \cite{Tokunaga}.

\begin{lem}[\cite{Tokunaga}]
A maximal outerplanar graph $G$ can be $4$-colored such that every
cycle of length $4$ in $G$ has all four colors.
\end{lem}

\begin{thm}
Let $G$ be a maximal outerplanar graph of order $n\geq 3$. If $t$ is
the number of vertices of degree $2$ in $G$, then $\gamma_{\times
2}(G)\leq \lfloor \frac{n+t}{2}\rfloor$.
\end{thm}

\bp Let $R=\{a_1, a_2, \cdots, a_t\}$ be the set of vertices of $G$
having degree $2$. For each $a_i$, let $b_i$ be one of its
neighbors. We construct a graph $G'$ from $G$ by adding $t$ new
vertices $u_1, u_2, \cdots, u_t$, and $2t$ new edges $a_1u_1,
b_1u_1, a_2u_2, b_2u_2, \cdots, a_tu_t, b_tu_t$. Note that $G'$ is
also a maximal outerplanar graph. By Lemma~3.4, $G'$ can be
$4$-colored such that every cycle of length $4$ in $G'$ has all four
colors. Let $V_p$ be the set of vertices assigned color $p$, where
$p=1, 2, 3, 4$. Choosing two suitable color classes, say $V_1$ and
$V_2$, such that $|V_1|+|V_2|\leq \lfloor \frac{n+t}{2}\rfloor$. Let
$D=V_1\cup V_2$.

Let $D\cap \{u_1, u_2, \cdots, u_t\}=\{u_1', u_2', \cdots, u_k'\}$,
where $k\leq t$. For any $u_i'$, at least one of its neighbors, say
$v_i$, does not belong to $D$. We construct a new set
$D'=(D\setminus \{u_1', u_2', \cdots, u_k'\})\cup \{v_1, v_2,
\cdots, v_k\}$.

Now, for any vertex $x$ of $G'$ of degree at least $3$. Let $x_1,
x_2$ and $x_3$ be consecutive vertices of $N_{G'}(x)$ in this order.
Clearly, $xx_1x_2x_3$ forms a cycle of length $4$, and two of $\{x,
x_1, x_2, x_3\}$ belong to $D$. Since there is at most one vertex of
degree two in $\{x, x_1, x_2, x_3\}$, at least two of those vertices
of degree at least three belong to $D'$. It means that $x$ is double
dominated by $D'$. And then, $D'$ is a double dominating set of $G$.
\ep

On the other hand, for a maximal outerplanar graph $G$ of order at
least four, it is easy to see that all vertices of $G$ of degree at
least three form a double dominating set of $G$. Hence we have the
following conclusion.

\begin{obs}
Let $G$ be a maximal outerplanar graph of order $n\geq 4$. If $t$ is
the number of vertices of degree $2$ in $G$, then $\gamma_{\times
2}(G)\leq n-t$.
\end{obs}

Clearly, each graph of $\mathscr{U}$ achieves equalities in
Theorem~3.5 and Observation~3.6. That is, the bounds of Theorem~3.5
and Observation~3.6 are also tight. The following conclusion is
immediate from Theorem~3.5 and Observation~3.6.

\begin{cor}
Let $G$ be a maximal outerplanar graph of order $n\geq 4$. If $t$ is
the number of vertices of degree $2$ in $G$, then
$$\gamma_{\times
2}(G)\leq \left\{
          \begin{array}{ll}
            \lfloor \frac{n+t}{2}\rfloor, & \hbox{if $t<\frac{n}{3}$;} \\
            n-t, & \hbox{otherwise}
          \end{array}
        \right.
        $$
\end{cor}

Next, we compare the bound in Corollary~3.3 with that in
Corollary~3.7. From Proposition~2.2, we know that there are at most
$\frac{n}{2}$ vertices of degree $2$ in a maximal outerplanar graph.
It is easy to see that $\lfloor \frac{n+t}{2}\rfloor<\lfloor
\frac{2n}{3}\rfloor$ when $t<\frac{n}{3}$, and $n-t<\lfloor
\frac{2n}{3}\rfloor$ when $\frac{n}{3}<t\leq \frac{n}{2}$. Thus, the
bound in Corollary~3.7 is better than that in Corollary~3.3.

Based on the above analysis, it is natural to consider the following
question: can the upper bound for Corollary~3.7 be improved for the
striped maximal outerplanar graph? From the above result, we know
that $\lfloor \frac{n+t}{2}\rfloor <n-t$ when $t$ is far less than
$n$. Then combining Proposition~2.1, we have the following.

\begin{cor}
Let $G$ be a striped maximal outerplanar graph with $n\geq 3$
vertices. Then, $\gamma_{\times 2}(G)\leq \lfloor
\frac{n}{2}\rfloor+1$.
\end{cor}

Next, we construct a striped maximal outerplanar graph $G$ of order
$n\geq 6$. Let $C=a_1a_2\cdots a_{q-1}a_qb_qb_{q-1}\cdots b_2b_1a_1$
is the unique Hamiltonian cycle of $G$, where $q\equiv 1(mod$ $4)$
or $q\equiv 3(mod$ $4)$. We know that $C$ is the boundary of the
outer face of $G$. Let $E_1=\{a_2b_2, a_3b_3, \cdots,
a_{q-1}b_{q-1}\}$, $E(G)=E(C)\cup E_1\cup \{a_1b_2, a_2b_3, a_4b_3,
a_5b_4, a_5b_6, a_6b_7, a_8b_7, a_9b_8, $\\$\cdots, a_{q-4}b_{q-3},
a_{q-3}b_{q-2}, a_{q-1}b_{q-2}, a_qb_{q-1}\}$ when $q\equiv 1(mod$
$4)$; and $E(G)=E(C)\cup E_1\cup \{a_1b_2, a_2b_3, a_4b_3, a_5b_4,
a_5b_6, a_6b_7, a_8b_7, a_9b_8, \cdots, a_{q-3}b_{q-4},
a_{q-2}b_{q-3}, a_{q-2}b_{q-1}, a_{q-1}b_q\}$ when \\$q\equiv 3(mod$
$4)$ (see Fig.2). Note that $|G|\equiv 2(mod$ $4)$. Let
$\mathscr{A}$ be a family consisting of all such graph $G$. We will
show that each graph of $\mathscr{A}$ attains the bound of
Corollary~3.8.

\begin{center}
  \includegraphics[width=6in]{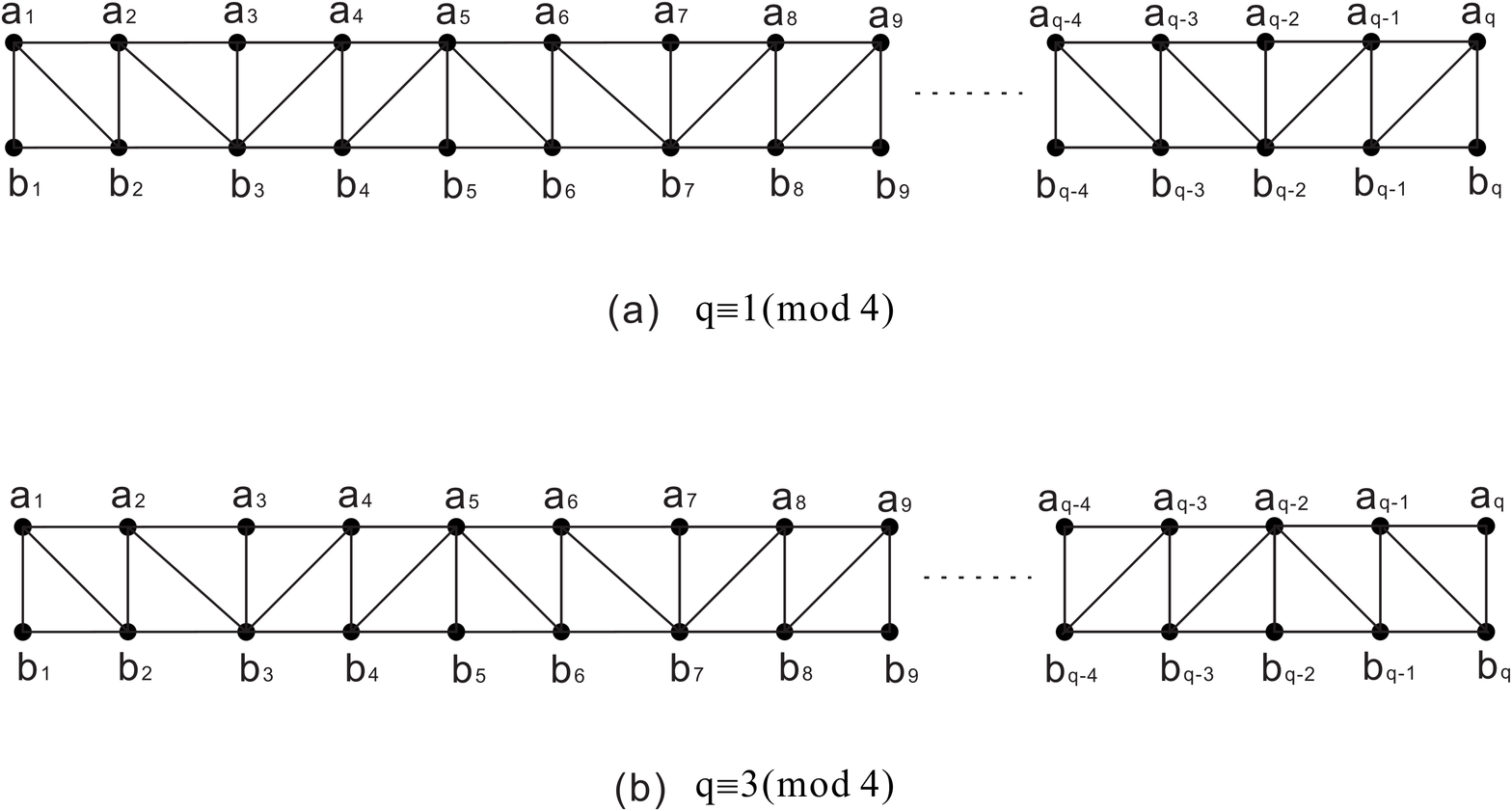}
 \end{center}
\qquad \qquad \qquad \qquad \qquad \qquad \qquad \qquad \qquad
{\small Fig.2} $\\$

\begin{center}
  \includegraphics[width=4in]{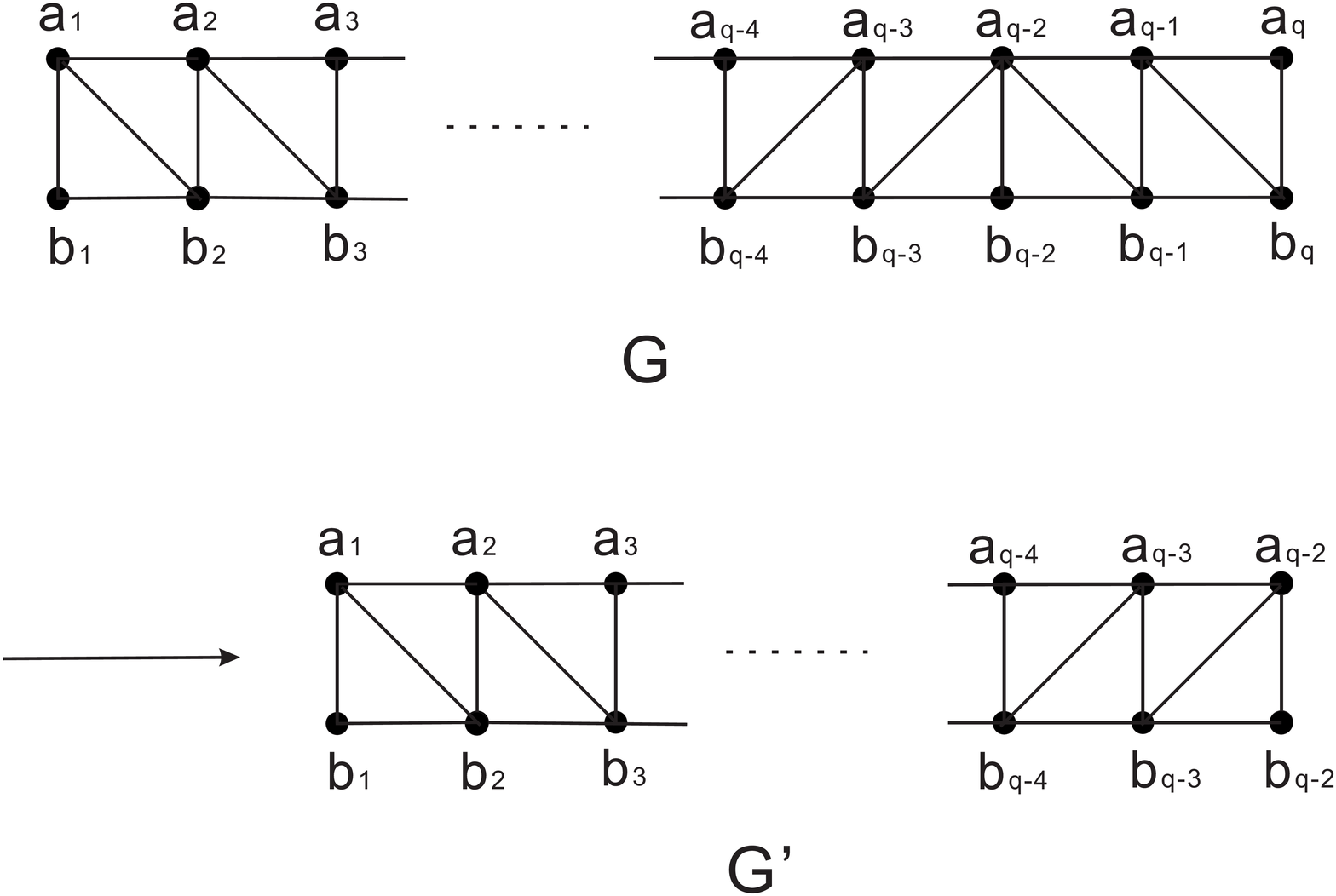}
 \end{center}
\qquad \qquad \qquad \qquad \qquad \qquad \qquad \qquad \qquad
{\small Fig.3} $\\$

\begin{thm}
If a striped maximal outerplanar graph $G$ belongs to $\mathscr{A}$,
then $\gamma_{\times 2}(G)=\frac{|G|}{2}+1$.
\end{thm}

\bp If $n=6$, the result holds. This establishes the base cases. So
we let $n\geq 10$ (Note that $n\equiv 2(mod$ $4)$). Assume that $G$
is a striped maximal outerplanar graph of $\mathscr{A}$ with minimum
order, such that $\gamma_{\times 2}(G)<\frac{|G|}{2}+1$.

Let $C=a_1a_2\cdots a_{q-1}a_qb_qb_{q-1}\cdots b_2b_1a_1$ be the
unique Hamiltonian cycle of $G$, and $S$ be a $\gamma_{\times
2}$-set of $G$. We know that $G$ has exactly two vertices of degree
two, without loss of generality, suppose that $a_q$ is one of them.
Clearly, $|N[a_q]\cap S|\geq 2$, and $S_1=(S\setminus N[a_q])\cup
\{a_{q-1}, b_q\}$ is still a $\gamma_{\times 2}$-set of $G$.
Moreover, $|N[b_{q-2}]\cap S_1|\geq 2$, and $S_2=(S_1\setminus
N[b_{q-2}])\cup \{a_{q-2}, b_{q-3}\}$ is also a $\gamma_{\times
2}$-set of $G$.

Let $G'=G-\{a_{q-1}, b_{q-1}, a_q, b_q\}$ (see Fig.3). Note that
$G'$ is also a striped maximal outerplanar graph which belongs to
$\mathscr{A}$, and $S_2\setminus \{a_{q-1}, b_q\}$ is a double
dominating set of $G'$. Thus, $\gamma_{\times 2}(G')\leq
\gamma_{\times
2}(G)-2<\frac{|G|}{2}-1=\frac{|G'|+4}{2}-1=\frac{|G'|}{2}+1$. It
contradicts the assumption that $\gamma_{\times
2}(G')=\frac{|G'|}{2}+1$. \ep

From the above result, we know that Corollary~3.8 is tight for
striped maximal outerplanar graphs.


\end{document}